\documentclass[12pt]{article}
\usepackage{graphicx}
\usepackage{amsfonts}
\newcommand{\pic}[4]{\vspace{1ex}\setlength{\unitlength}{1cm}
\begin{picture}(0,#3)(5,.5)
\put(#2){\includegraphics[#4]{#1.eps}}
\end{picture}\vspace{1ex}}

\newtheorem{prop}{Proposition}[section]

\newtheorem{theo}[prop]{Theorem}
\newtheorem{lemm}[prop]{Lemma}
\newtheorem{conj}[prop]{Conjecture}
\newtheorem{defi}[prop]{Definition}

\def\proof{\par\vspace{-1ex}\noindent\textit{Proof. }}
 
\def\C{\mathbb{C}}
\def\cg{\overline} 

\def\D{\mathbb{D}}

\def\G{\mathcal{G}}
\def\R{\mathbb{R}}
\def\S{{\mathcal{S}}}

\def\qed{\hskip1em\raise3.5pt\hbox{\framebox[2mm]{\ }}}

\def\im{\,{\rm Im}\,}

\begin{document}

\begin{center} {\LARGE
    Gears, Pregears and Related Domains}
\end{center}
\begin{center} 
  \today
\end{center}
\begin{center}
    Philip R.  Brown\footnote{Partially supported by CONACyT grant
      166183} \\
    R. Michael Porter\footnotemark[1] 
\end{center}
 
\noindent Abstract.  We study conformal mappings from the unit disk to
one-toothed gear-shaped planar domains from the point of view of the
Schwarzian derivative. Gear-shaped (or ``gearlike'') domains fit into
a more general category of domains we call ``pregears'' (images of
gears under M\"obius transformations), which aid in the study of the
conformal mappings for gears and which we also describe in
detail. Such domains being bounded by arcs of circles, the Schwarzian
derivative of the Riemann mapping is known to be a rational function
of a specific form. One accessory parameter of these mappings is
naturally related to the conformal modulus of the gear (or pregear)
and we prove several qualitative results relating it to the principal
remaining accessory parameter. The corresponding region of univalence
(parameters for which the rational function is the Schwarzian
derivative of a conformal mapping) is determined precisely.


\medskip

\section{Introduction \label{sec:intro}}
A special case of a circular quadrilateral is a \textit{gear domain}
with one tooth: a starlike open set in the complex plane bounded by
arcs of two circles centered at the origin and segments of two lines
passing through the origin. A related family of domains which we call
\textit{pregear domains}, are those which are M\"obius transformations
of gear domains.

In \cite{Goo}, the Riemann mapping of the unit disk onto a gear
domain, fixing the origin, is expressed as the solution of a first
order linear differential equation (see (\ref{eq:goodman}) below) which is
derived by making use of the starlike property in the unit disk and
the boundary behavior of the mapping.

An alternative approach, which we follow here, follows the more
general construction of conformal maps onto circular polygons,
expressing the mapping as a solution of the equation which prescribes its
Schwarzian derivative, as expounded in \cite{BG,Br2,BrP}, \cite[p.\
70]{DT}, \cite{Hi,KP2}, \cite[p.\ 198]{Neh}, and \cite{Phoro} for
circular quadrilaterals. In the present case of gear domains with one
tooth, after normalizing the location of the prevertices, the
Schwarzian derivative $R_{t,\lambda}(z)$ (see
(\ref{eq:Rtlambda})--(\ref{eq:psi0tpsi1t}) below) contains the two
unknown parameters $t$ (which determines the prevertices) and
$\lambda$ (an auxiliary parameter), which in turn determine via the
conformal mapping the two natural geometric quantities which specify
the gear domain (Figure \ref{fig:onetooth}); namely, the ratio of the
radii of the outermost and  the innermost circle,
which we call the \emph{gear ratio} $\beta$, and half the angle
between the rays that form the straight boundary segments, which we
call the \emph{gear angle} $\gamma$.

We describe one-tooth gear domains, summarize the relevant theory from
\cite{Goo}, and work out the formula for $R_{t,\lambda}(z)$ in
Section~\ref{sec:quads}.  The functional relationship between the
parameters $t$ and $\lambda$ and the corresponding gear domain is
analyzed Section~\ref{modules}, where for convenience we work with the
conformal module $M(t)$ of the gear domain. We prove that when $t$ is
fixed, for each $\gamma$ there is at most one value of $\lambda$ for
which a solution of (\ref{eq:SfRtlambda0}) is a gear. We also prove
that when $\gamma$ is fixed, $\beta$ is a monotonic function of $M(t)$
in the full range $0<t<\pi/2$. On the other hand, when $\beta$ is
fixed, $M(t)\to 0$ as $\gamma\to 0$ or $\gamma \to\pi$, which leads to
our conjecture that there are exactly two gears (corresponding to two
different values of $\gamma$) with module $M(t)$, provided $t$ is
below a threshold value $t_\beta$; for $t=t_\beta$ there is only one
gear with module $M(t)$ and when $t$ is above the threshold there are
no gears at all with this module. 

A basic property of the Schwarzian derivative is its invariance under
M\"obius transformations. Consequently, the solutions $f$ obtained by
solving (\ref{eq:sdeq}) or (\ref{eq:SfRtlambda0}) with the classical
normalizations \mbox{$f(0)\!=\!0$} and \mbox{$f'(0)\!=\!1$} are, in
general, M\"obius transforms of one-tooth gear domains. In
Section~\ref{sec:gearregion} we give a full geometric description of
these ``pregear'' domains and regard them to be of independent
interest from the point of view of conformal mapping.  The family of
pregear domains is seen to be bounded by ``degenerate pregear domains'', which
are M\"obius transformations of certain unbounded rectilinear
quadrilaterals. From this observation we are able to determine
precisely the boundary of the region in the $(t,\lambda)$-plane for
which $R_{t,\lambda}$ is the Schwarzian derivative of a univalent
function (i.e., a conformal mapping to a gear) .

\section{Schwarzian derivative and accessory\\ parameters \label{sec:quads}}
  
\subsection{Gear domains}

In this study a \textit{gearlike domain} (or gear domain) is
a starlike open set $G$ in the complex plane $\C$ bounded by   
arcs of circles centered at the origin and segments of lines passing
through the origin. Occasionally for reasons of normalization of mappings
we may use the same term for a translate of a gear domain, in which case
we will clarify that the ''gear center'' may not be the origin.
 The fundamental study of gearlike domains was
initiated in \cite{Goo} and further results have appeared in 
\cite{BPea,Br3,Ni,Pea}.
 
\begin{figure}\centering
\pic{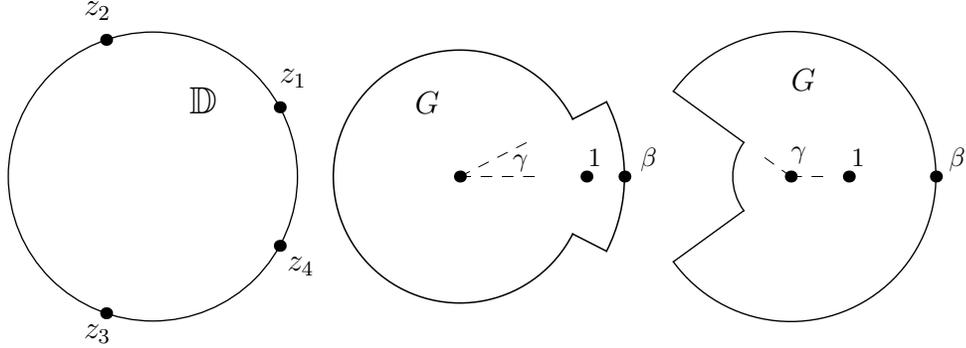}{-1,0}{4}{}
\caption{Gear parameters.  }\label{fig:onetooth}
\end{figure}

We will consider in particular bounded gearlike domains $G$ with a
single ``tooth'' as in Figure \ref{fig:onetooth}.  The interior angles
are defined by $\pi\alpha_i$ for $i=1,2,3,4$, where
$\alpha_1=\alpha_4=1/2$ and $\alpha_2=\alpha_3=3/2$. We will assume that
$G$ is symmetric in the real axis, and that the corresponding
prevertices of the conformal mapping are of the form
\[ z_1=e^{it_1},\  z_2=e^{it_2},\  z_3=e^{-it_2},\  z_4=e^{-it_1},
\]
$0<t_1<t_2<\pi$, as can always be obtained by a preliminary
transformation of $\D$.  The straight edges of $G$ will be be referred
to as the \emph{tooth edges}, which when prolonged meet at the
\emph{gear center} at an angle $2\gamma$, where $\gamma$ will be
termed the \emph{gear angle}.  The two edges of $\partial G$ which are
not tooth edges are arcs of circles centered at the gear center: the
\emph{A-arc} ending at angles of $3\pi/2$ and the \emph{B-arc} ending
at angles of $\pi/2$. The quotient $\beta$ of the radius of the B-arc
to that of the A-arc is the \emph{gear ratio} of $G$. We write
$G_{\beta,\gamma}$ for the standard gear domain with gear parameters
$\beta,\gamma$ which is centered at the origin and has A-arc of radius
$1$.
  
%

Most of the research which has been done on gearlike domains is based on
results in \cite{Goo}.  We state here the following particular case.
\begin{prop}\label{prop:goodman}
  A necessary and sufficient condition for $f(z)$ to be a univalent
  mapping of the unit disk onto a one-tooth gear domain, satisfying
  the normalizations $f(0)=0$ and
\begin{equation}\label{derivatives}
f''(0)=2f'(0)(\cos t_2-\cos t_1),
\end{equation}
where the prevertices $e^{\pm it_1}$ map to the vertices with interior
angles $\pi/2$ and the prevertices $e^{\pm it_2}$ map to the vertices
with interior angles $3\pi/2$, is that
\begin{equation}\label{eq:goodman}
 f'(z)= \frac{1}{z}\left(\frac{z^2-2z\cos t_1+1}{z^2-2z\cos t_2+1}\right)^{1/2} f(z).
\end{equation}
Furthermore, the gear ratio $\beta$ and the gear angle $\gamma$ are
determined by the following two integral formulas:
\[ \log \beta = \int_{t_1}^{t_2}\sqrt{\frac{\cos\theta-\cos t_2}
                                {\cos t_1 -\cos\theta  }}\,d\theta, \quad
\gamma = \int_0^{t_1}\sqrt{\frac{\cos\theta-\cos t_2}
                                {\cos\theta  -\cos t_1 }  }\,d\theta.
\]
\end{prop}

\proof The necessity of (\ref{eq:goodman}) is a consequence of
\cite[Theorem 2]{Goo}. Goodman also derived the formulas for $\beta$
and $\gamma$ and gave formulas for the coefficients $\{b_n\}_{n\ge2}$
of the Maclaurin series of $f(z)$ in terms of the coefficient
$b_1$. In particular, he obtained (\ref{derivatives}). The sufficiency
of (\ref{eq:goodman}) is a consequence of the geometry and symmetry of
a one-tooth gear domain: divide the unit circle $\partial\D$ into four
arcs separated by the points $z_1,z_2,z_3,z_4$. It follows from
(\ref{eq:goodman}) as in the proofs of \cite[Lemmas 1 and 2]{Goo} that
as $\theta$ increases, $|f|$ is constant and $\arg f$ is increasing
for $-t_1<\theta<t_1$ and $t_2<\theta<2\pi-t_2$; similarly $\arg f$ is
constant and $|f|$ is decreasing for $t_1<\theta<t_2$, while $\arg f$
is constant and $|f|$ is increasing for $-t_2<\theta<-t_1$.  As a
result of the singularities in (\ref{eq:goodman}), the interior angles
of $f(\partial\D)$ at $f(e^{\pm it_1})$ must be $\pi/2$ and the
interior angles at $f(e^{\pm it_2})$ must be $3\pi/2$.  Now if $f$
were not univalent, then an examination of a few possibilities for
which the properties above are satisfied would show that the winding
number of $f(\partial\D)$ about the origin must be greater than 1,
which means that $f$ could not have a simple zero at the
origin. However, this contradicts the property
$\lim_{z\to0}z\,f'(z)/f(z)=1$.  \qed
 
While it is straightforward to solve (\ref{eq:goodman}) numerically,
from these relations it is difficult to find the values of $t_1,t_2$
corresponding to a pair $\beta,\gamma$.

\subsection{The Schwarzian derivative of a gear mapping}

Although in principle the operator $zf'/f$ considered by Goodman is
simpler than the Schwarzian derivative
$\S_f=(f''/f')'-(1/2)(f''/f')^2$, we will base our study of gear
domains on the latter to take advantage of the rich theory which has
developed around it, as referred to in the introduction.
 
We require some general facts about
conformal mappings of circle polygon domains.
Let $\D=\{|z<1\}$ denote the unit disk.  For a general circle polygon
domain $D$ with interior angles $\pi\alpha_k$ at the vertices
$w_k$, write $a_k=(1-\alpha_k^2)/2$. Let $f\colon\D\to D$
be a conformal mapping.  Then the Schwarzian derivative $\S_f$  is a
rational function of the form 
\begin{equation} \label{eq:schgen}
     \S_f(z) = z^{-2}\sum_{k=1}^n \bigg( 
 \frac{a_k z_k z}{(z-z_k)^2} +  i r_k\frac{z+z_k}{z-z_k}   \bigg),
\end{equation}
where $z_k=f^{-1}(w_k)\in\partial\D$ ($1\le k\le n$) are the prevertices of $D$ with respect
to the mapping $f$, and $r_k\in\R$ are additional
accessory parameters which satisfy the relations
\begin{equation} \label{eq:rcondgen} \sum_{k=1}^n r_k = 0, \quad
                 \sum_{k=1}^n z_k(a_k+2ir_k) = 0,
\end{equation} 
These relations assure that $f$ sends the boundary $\partial\D$ to a
union of circular arcs, and further imply that the singularity of
(\ref{eq:schgen}) at the origin is removable. See \cite{Phoro} for a proof.
For one-tooth gear domains this specializes as follows. 

\begin{theo} \label{th:ratfuncs} Let $G$ be a one-tooth gear domain and let
  $f\colon\D\to G$ be a conformal mapping.  Suppose that $f$ is
  symmetric in the real axis.  Then there are unique values
  $t_1,t_2,\lambda$ ($0<t_1<t_2<\pi$, $\lambda\in\R$) such that Schwarzian derivative
  $\S_f$ of $f$ can be expressed as 
\begin{equation}    \label{eq:sdeq}
S_f=R_{t_1,t_2,\lambda}
\end{equation}
 where
\begin{equation}    \label{eq:ratfuncs}
 \frac{1}{2} R_{t_1,t_2,\lambda}(z)=  \psi_{0,(t_1,t_2)}(z) - \lambda \psi_{1,(t_1,t_2)}(z) 
\end{equation}
with
\begin{equation}\label{eq:psi1}
   \psi_{1,(t_1,t_2)}(z)  = \frac{4 (\cos t_2-\cos t_1) }
{\left(z^2-(2\cos t_1) z  +1\right)  \left(z^2-(2 \cos t_2) z +1\right)}
\end{equation}
and
\begin{equation}\label{eq:psi0}
  \psi_{0,(t_1,t_2)}(z)  =  \frac{ c_{40} z^4  +  c_{30}z^3 +  c_{20} z^2 + c_{10}z  +   c_{00} }
{\left(z^2-(2\cos t_1) z  +1\right)^2  \left(z^2-(2 \cos t_2) z +1\right)^2},
\end{equation}
\begin{eqnarray*}
 c_{00} &=&  c_{40} \ =\ \frac{ 3 \cos 2t_1-5 \cos 2 t_2+2}{8},   \\
 c_{10} &=&   c_{30} \ =\ 3 \sin ^2t_1 \cos t_2-5 \cos t_1 \sin ^2t_2 ,   \\[1ex]
 c_{20} &=& \frac{(\cos 2 t_1) (11-2 \cos 2 t_2)-13 \cos 2 t_2+4}{4} . 
\end{eqnarray*}
\end{theo}

\proof Via the assumed symmetry $f(z)=\cg{f(\cg{z})}$ we have
$\S_f(z)=\cg{\S_f(\cg{z})}$.  From this and the general form
(\ref{eq:schgen}) with $n=4$ it follows (using
$z_1=\exp(it_1)=\cg{z_4}$, $z_2=\exp(it_2)=\cg{z_3}$) that $r_1=-r_4$, $r_2=-r_3$,
so with the first relation of (\ref{eq:rcondgen}) the general formula
reduces to
\begin{eqnarray*} z^2 \S_f(z)=
 && \frac{3}{8}\left(
   \frac{ z {z_1}}{(z-{z_1})^2}+\frac{ z\cg{z_1}}{(z-\cg{z_1})^2} \right)
   -\frac{5}{8}\left(
  \frac{ z {z_2}}{(z-{z_2})^2}+\frac{ z \cg{z_2}}{(z-\cg{z_2})^2} \right) \\
&& +i {r_1} \left(\frac{z+{z_1}}{z-{z_1}}-\frac{z+\cg{z_1}}{z-\cg{z_1}}\right)+i {r_2}
   \left(\frac{z+{z_2}}{z-{z_2}}-\frac{z+\cg{z_2}}{z-\cg{z_2}}\right)
\end{eqnarray*} 
 We introduce
a new real parameter $\lambda$, determined by
 \begin{equation}
  r_1 = \frac{ \lambda + (3/8)\cos t_1}{2\sin t_1}, \quad
  r_2= \frac{ \lambda - (5/8)\cos t_1}{2\sin t_2} ,
\end{equation}
and the second relation of (\ref{eq:rcondgen}). After
substituting $\exp\pm it_1$ and $\exp\pm it_2$ for the prevertices of
$f$, the result follows by algebraic manipulation (see
\cite[Appendix]{Phoro} for further details). Clearly $t_1$,
$t_2$ are determined by the given prevertices of $f$, and the uniqueness
of $\lambda$ follows from the explicit formula. \qed

\noindent\emph{Symmetrization of prevertices.} We are interested in the combinations of
accessory parameters $(t_1,t_2,\lambda)$ for which $f$ is a 
gear mapping. We have the flexibilty to relocate the prevertices $t_1$
and $t_2$, as we explain in the following paragraphs.

The \textit{pullback} of a gear mapping $f$ via the self-mapping
\begin{equation}  \label{eq:Tq}
  T_q(z) = \frac{z-q}{-qz+1}.
\end{equation}  of $\D$ is by definition the Schwarzian
derivative of the function $f\circ T_q$, also defined in $\D$.  Since
the composition is also a gear mapping, by Theorem \ref{th:ratfuncs}
and the Chain Rule we obtain the form of the  pullback of $\S_f$,
\begin{equation}\label{eq:pullback} 
  R_{t_1,t_2,\lambda}(z) =
     R_{t_1^*,t_2^*,\lambda^*}(z^*)\, T_q'(z)^2,
\end{equation}
for some $\lambda^*\in\R$, where $\S_f=R_{t_1^*,t_2^*,\lambda^*}(z^*)$,
$e^{it_1^*}=T_q(e^{it_1})$, $e^{it_2^*}=T_q(e^{it_2})$, $z^*=T_q(z)$.
It is not difficult to see that the ``$\psi_0$'' and ``$\psi_1$''
parts of the Schwarzian derivative pull back independently because
of the differing degrees of the polynomials in their denominators:
\begin{eqnarray} 
 \psi_{0;t_1,t_2}(z) &=& \psi_{0;t_1^*,t_2^*}(z^*)\, T_q'(z)^2, \nonumber\\
 \psi_{1;t_1,t_2}(z) &=& \psi_{1;t_1^*,t_2^*}(z^*)\, T_q'(z)^2.  
\end{eqnarray}
In fact, the following holds.

\begin{prop} \label{prop:pullback} 
  The pullback of the Schwarzian derivative $\S_f$ of a gear mapping
  respects the auxiliary parameters in the sense that
  $\lambda^*=\lambda$ in (\ref{eq:pullback}).
\end{prop}

\proof This is an explicit but tedious calculation making use of the
elementary formula
\[ T_q(a)-T_q(b) = T_q'(a)^{1/2}T_q'(b)^{1/2}(a-b) =
\frac{(1-q^2)(a-b)}{(1-qa)(1-qb)}
\]
which can be derived immediately from the the formula for the
derivative $T_q'(z)=(1-q^2)/(-qz+1)^2$ (and choosing say
$\sqrt{1-q^2}/(-qz+1)$ as the indicated square root).  \qed

The following result shows that by means of an appropriate $T_q$ it is
possible to make the prevertices symmetric in the imaginary axis. This
will be useful for our considerations of modules in Section
\ref{modules}.

\begin{prop} \label{prop:Tq} 
Let $0<t_1<t_2<\pi$,  $z_1=e^{it_1}$, $z_2=e^{it_2}$. Define   
\[  q = \frac{z_1+z_2-2\sqrt{\im z_1\,\im z_2}\sqrt{z_1z_2}}{1+z_1z_2} 
\]
 where  $\sqrt{z_1z_2}=e^{i(t_1+t_2)/2}$. Then there is  a unique $0<t<\pi/2$ such that
\[   T_q(z_1)=e^{it}, \quad   T_q(z_2)=e^{i(\pi-t)} .
 \] 
\end{prop}

\proof First we note that $\cg{q}=q$ because of the relations
$z_1\cg{z_1}=z_2\cg{z_2}=1$.  It follows from the definition (\ref{eq:Tq}) that
$|T_q(z_1)|=|T_q(z_2)|=1$, and then we calculate that
 \begin{eqnarray*}
 T_q(z_1) &=&  \frac{z_1-z_2+2\sqrt{\im z_1\,\im z_2}\sqrt{z_1z_2}}{1-z_1z_2}, \\
 T_q(z_2) &=&  \frac{-z_2-z_1+2\sqrt{\im z_1\,\im z_2}\sqrt{z_1z_2}}{1-z_1z_2}. 
\end{eqnarray*}
Using the same reasoning we used to show  $\cg{q}=q$ it follows that $T_q(z_2)=-\cg{T_q(z_1)}$,
so we may take $t=\arg T_q(z_1)$. 

 Since $T_q$ fixes $\pm 1$ and thus conserves the 
order of $1,e^{it_1},e^{it_2},-1$ along $\partial\D$, it follows that $T_1'(0)>0$, which implies
that $|q|<1$.
\qed

By Proposition~\ref{prop:pullback} we can always
express $ R_{t_1,t_2,\lambda}$ as
\begin{equation}  \label{eq:Rt1t2lambda}
    R_{t_1,t_2,\lambda} = (R_{t,\pi/2-t,\lambda}\circ T_q) (T_q')^2.
\end{equation} 
for a parameter $t$ uniquely determined by $t_1,t_2$ according to
Proposition~\ref{prop:Tq}.  We thus reduce the study of gear mappings to the special case
\begin{equation} \label{eq:RtlambdafromRt1t2}
 R_{t ,\lambda} := R_{t,\pi/2-t,\lambda}.
\end{equation}
Explicitly, we have a much simpler Schwarzian derivative
\begin{eqnarray}  \label{eq:Rtlambda}
    \frac{1}{2}R_{t ,\lambda}(z) &=& \psi_{0,t}(z) - \lambda
    \psi_{1,t}(z)
\end{eqnarray}
where
 \begin{eqnarray}  \label{eq:psi0tpsi1t} 
  \psi_{0,t}(z) 
   &=&    \frac{ (\sin^2t)(z^4-(16\cos t)z^3 + (4+2\cos2t) z^2 - (16\cos t) z + 1) }
           {2(z^4-(2 \cos2t)z^2+1)^2}, \nonumber \\
  \psi_{1,t}(z)  &=&    \frac{-8\cos t}{z^4-(2\cos2t)z^2+1}.  \label{eq:psit}
\end{eqnarray}
We make the precautionary observation that even though the prevertices
are now symmetric in both coordinate axes, the Schwarzian derivative
$R_{t ,\lambda}$ is not symmetric in the imaginary axis.

\section{Conformal modules of gear domains}\label{modules}


We discuss some of the relations among $t,\lambda,\beta,\gamma$.
First we will treat $t$ as fixed, for the following reason.  The
\textit{conformal module} $M(G_{\beta,\gamma})>0$ of any gear (or
pregear) with prevertices $\pm e^{\pm it}$ is by definition the
conformal module of the unique conformally equivalent rectangle
$(0,1,1+\tau,\tau)$ with imaginary $\tau$; i.e., $\tau=iM(t)$.  Thus
we can write $M(t)=M(G_{\beta,\gamma})$.  Note that
\begin{equation}\label{eq:limitM(t)}
 \lim_{t\to0}M(t)=0; \quad  \lim_{t\to\pi/2}M(t)=\infty .
\end{equation}


\begin{defi}\rm
We will say that the rational function $R_{t,\lambda}$ is
\textit{gearlike} if there is a solution $f$ of
\begin{equation}    \label{eq:SfRtlambda0}
   \S_f = R_{t,\lambda}
\end{equation}
that is a univalent mapping onto a gear.
\end{defi}

\begin{lemm} \label{lemm:lambdafromgamma}
  Fix $t$.  Then for each $\gamma$ there is at most one value of
  $\lambda$ for which $R_{t,\lambda}$ is
  gearlike. 
\end{lemm}
\proof Different values of $\lambda$ for which $R_{t,\lambda}$ is
gearlike would correspond to different values of $\beta$. If two
gears have the same gear angle $\gamma$ but different values of $\beta$, one must 
be contained within the other. By   the
well-known monotonicity of conformal modules of topological
quadrilaterals \cite{LV} it is not possible for these
gears to have the same $M(t)$.  \qed

\begin{prop} \label{prop:limitmodules}
 (i) Let $0<\gamma<\pi$. Then 
 $M(G_{\beta,\gamma})\to\infty$ as $\beta\to1$, while
 $M(G_{\beta,\gamma})\to0$ as $\beta\to\infty$.
(ii) Let $\beta>1$. Then 
 $M(G_{\beta,\gamma})\to0$ as $\gamma\to0$ or $\gamma\to\pi$.
\end{prop} 

\proof (i) The Euclidean separation of the vertical sides of the
topological quadrilateral $G_{\beta,\gamma}$ tends to 0 as $\beta\to1$
while the horizontal sides are bounded away from one another, hence
the conformal module tends to $\infty$ (\cite[Lemma 4.1]{LV}).
Similarly, as $\beta\to\infty$, one may rescale the quadrilateral to
see that the separation of the horizontal edges
tends to 0, so the conformal module tends to $0$.
 
(ii) This limit is harder to see because although the horizontal edges
in $G_{\beta,\gamma}$ become arbitrarily close, the rightmost vertical
edge also degenerates (the Riemann map with $0\to0$ takes all four
vertices close to 1).  We can decompose the gear as
$G_{\beta,\gamma}=\D \cup \{e^{i\theta}\colon\ |\theta|<\gamma\}\cup
D_1$ where $D_1$ is the tooth. (Here $\D$ is not a quadrilateral, just
a bigon attached to a vertical side of $D_1$ but it can be thought of
as a limiting case of a quadrilateral.) By monotonicity of conformal
modulus,
\[  \frac{1}{M(G_{\beta,\gamma})} \ge \frac{1}{M(D_1)}.
\]
Since $M(D_1)\to0$ as $\gamma\to0$, also $M(G_{\beta,\gamma})\to0$.

Now let $\gamma$ tend to $\pi$. The arc of $\partial\D$ from
$\gamma$ to $2\pi-\gamma$ is a vertical edge which is very small,
since there are curves of length less than $2(\pi-\gamma)$ interior to
$G_{\beta,\gamma}$ joining the vertical edges (which are the tooth
edges, having length $\beta-1$) to each other. On the other hand, the
Euclidean area of $G_{\beta,\gamma}$ is more than $\pi$, the area of $\D$,
independently of $\gamma$.  By \cite[sec.\ 4.3]{LV},
\[ M(G_{\beta,\gamma}) < \frac{(2(\pi-\gamma))^2}{\pi} \to 0.
\]
as $\gamma\to\pi$.  \qed

As in Lemma \ref{lemm:lambdafromgamma}, the function $\beta\mapsto
M(G_{\beta,\gamma})$ is monotone since the gear grows with
$\beta$. By Proposition \ref{prop:limitmodules} we have the following.

\begin{prop}\label{prop:betafromt}   Let $\gamma\in(0,\pi)$.
Then for any\/ $t\in(0,\pi/2)$ there is a unique $\beta>1$ such
that the gear $G_{\beta,\gamma}$ has conformal module $M(t)$.  
\end{prop} 
 
Given $\beta$, the existence of a value $\gamma$ for which $M(G(\beta,\gamma)))$ is
maximal follows likewise from Proposition \ref{prop:limitmodules}.  The
affirmation that the conformal module is monotone for $\gamma$ above
and below this value is so far confirmed only by numerical evidence
which we will present in our forthcoming article \cite{BrP3}.

\begin{conj} \label{conj:betafromgamma}
  Let $\beta>1$. Then there is a value $t_\beta\in(0,\pi/2)$ such that
  for each $t\in(0,t_\beta)$ there are exactly two values of $\gamma$
  such that $M(G(\beta,\gamma))=M(t)$. (For $t=t_\beta$ there is
  exactly one, while for $t>t_\beta$ there are none.)
\end{conj}  
  
We now have in hand enough information about gears to prove
one of our main results.

\begin{theo} \label{theo:lambdainterval} For each $t$ in
  $(0,\pi/2)$ there are constants $\lambda_t^-$,
  $\lambda_t^+$ such that $R_{t,\lambda}$ is the Schwarzian derivative
  of a conformal mapping from $\D$ to a gear domain if and only if
  $\lambda_t^-<\lambda<\lambda_t^+$.
\end{theo}
\proof By Proposition \ref{prop:limitmodules} there exist
$\beta_0,\gamma_0$ such that $M(G_{\beta_0,\gamma_0})=M(t)$.  By the
Riemann Mapping Theorem and Theorem \ref{th:ratfuncs}, there exists
$\lambda_0\in\R$ such that $R_{t,\lambda_0}$ is the Schwarzian
derivative of a conformal mapping from $\D$ to $G_{\beta_0,\gamma_0}$.
Since the set of $\lambda\in\R$ for which $R_{t,\lambda}$ is the
Schwarzian derivative of a conformal mapping from $\D$ to a gear
domain is clearly open, we may consider the maximal open interval
$I_t=(\lambda_t^-,\lambda_t^+)$ which contains $\lambda_0$ and is
contained in this set.  Suppose that there is a sequence
$\lambda_n\to\lambda_t^+$ in $I_t$ for which the gear ratio $\gamma_n$
converges to some $\gamma_0\in(0,\pi)$.  The corresponding conformal
mappings $f_n\colon\D\to G_{\beta_n,\gamma_n}$ with
$\S_{f_n}=R_{t,\lambda_n}$, symmetric in $\R$, and with $f_n'(0)>0$,
all cover $\D$ and thus (perhaps on a subsequence) converge to a
mapping $f$ with $\D\subseteq f(\D)$.
Since $\lambda_t^+$ is a boundary parameter, we must have
$\beta_n\to1$ or $\beta_n\to\infty$ on a subsequence.  However, by
Proposition \ref{prop:limitmodules} this implies
$M(G_{\beta_n,\gamma_n})\to0$ or $M(G_{\beta_n,\gamma_n})\to\infty$,
which is absurd since the conformal module $M(t)$ is fixed. It follows
that $\gamma$ accumulates only to $0$ or $\pi$ as
$\lambda\to\lambda_t^+$, and the same holds as
$\lambda\to\lambda_t^-$.

However, the function $\lambda\mapsto\gamma$ (for the fixed $t$ under
consideration) is strictly monotone.  Indeed, suppose that the
Schwarzian derivatives $R_{t,\lambda_1}$ and $R_{t,\lambda_2}$
determined gears $G_{\beta_1,\gamma}$ and $G_{\beta_2,\gamma}$ with
the same $\gamma$ but $\beta_1<\beta_2$. Then
$G_{\beta_1,\gamma}\subseteq G_{\beta_2,\gamma}$ and by monotonicity
of conformal modules, $M(G_{\beta_2,\gamma})<M(G_{\beta_1,\gamma})$,
again contradicting the fact that both modules are equal to $M(t)$.
We thus see that $\gamma\to\pi$ as $\lambda\to\lambda_t^-$ and
$\gamma\to0$ as $\lambda\to\lambda_t^+$, or vice
versa. (Numerically we will see in \cite{BrP3} that $\lambda\mapsto\gamma$ is
actually decreasing.)
 
Now we can show that $R_{t,\lambda}$ is not the Schwarzian derivative
of a gear mapping for any $\lambda$ outside of the interval $I_t$.  Suppose indeed
that $\lambda_1$ were such a parameter.  As was the case for
$\lambda_0$, there is a maximal interval $I$ containing $\lambda_1$
and in which every $R_{t,\lambda}$ gives a gear mapping. As before,
the limits of the $\gamma$ values at the endpoints of $I$ must be $0$
and $\pi$, and these values thus range over the whole interval
$(0,\pi)$.  Therefore there exists a value $\lambda\in I$ such that
$R_{t,\lambda}$ produces a gear $G_{\beta,\gamma_0}$, where
$R_{t,\lambda_0}$ produces $G_{\beta_0,\gamma_0}$. Again by
monotonicity of the module we have that $\beta=\beta_0$.  Then
by the uniqueness of conformal mappings, $R_{t,\lambda}=R_{t,\lambda_0}$
and finally $\lambda=\lambda_0$. This proves the statement.
\qed

(The limiting rational expressions $R_{t,\lambda^\pm_t}$ correspond to
maps to degenerate gears, as will be discussed now).

\section{The region of gearlikeness\label{sec:gearlikeness}}\label{sec:gearregion}

\subsection{Pregear domains}\label{sec:pregears}

\begin{defi}\rm
  We will say that a domain $D$ is a \emph{pregear} when it is the
  image $D=T(G)$ of a one-tooth gear domain $G$ under a M\"obius
  transformation $T$.
\end{defi}

\begin{figure}[h!tb]\centering
\pic{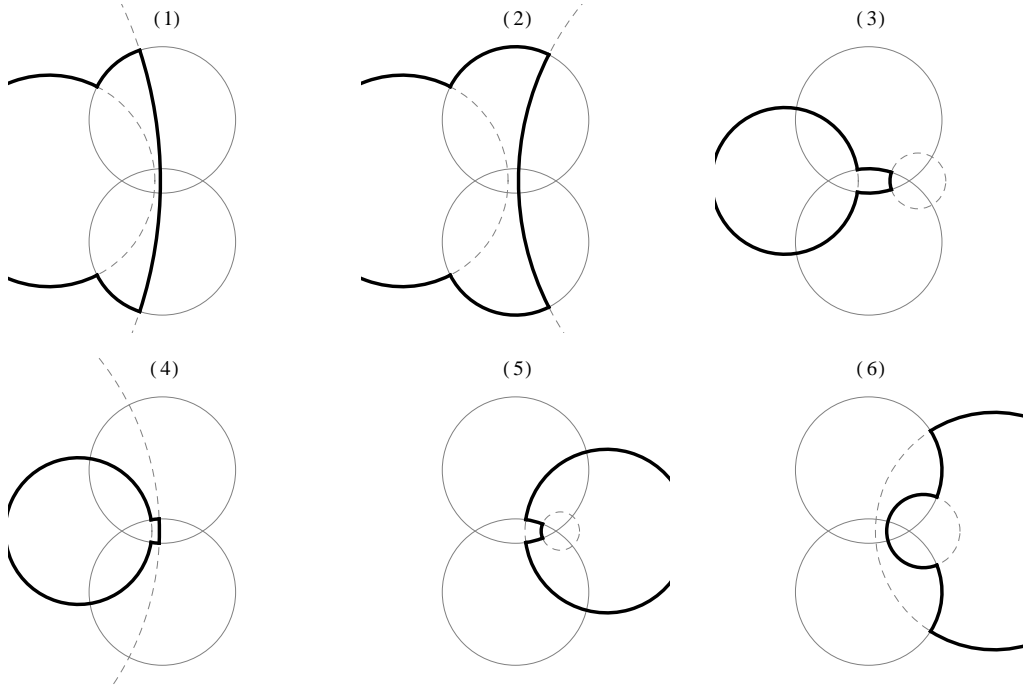}{-2,0}{8.5}{scale=.65} 
\caption{Examples of pregears. The circles $C^\pm$ containing the
  tooth edges are solid gray; the circles containing the A- and B-arcs
  in dotted gray. }
\label{fig:pregears1}
\end{figure}

The tooth edges as well as the A- and B-arcs of any pregear may be uniquely
identified by the corresponding interior angles.  Normally we will restrict the
discussion to pregears with the symmetries of the following
proposition.


\begin{prop} \label{prop:pregearcondition} Let $D$ be a circular
  quadrilateral which is symmetric in $\R$, has no vertices on $\R$,
  and has two interior angles equal to $\pi/2$ and two interior angles
  equal to $3\pi/2$. Assume that one tooth edge of $D$ lies in the
  upper and the other in the lower half-plane.  Then $D$ is a pregear
  if and only if the full circles $C^+$, $C^-$ containing the tooth edges
  intersect in two points.
\end{prop}
\proof Suppose that $C^+$ and $C^-$ intersect at the points $b^-$ and
$b^+$, which by symmetry are necessarily in $\R$. We can suppose that
$b^-<b^+$. If the A-arc passes between $b^-$ and $b^+$ as in diagrams
(5) and (6) in Figure \ref{fig:pregears1}, then $b^+$ is interior to
$D$. Otherwise, if the B-arc passes between $b^-$ and $b^+$ as in
diagrams (1), (2), (3) and (4) in Figure \ref{fig:pregears1}, then
$b^-$ is interior to $D$. Let $T$ be a M\"obius transformation such
that $T^{-1}(z)=(z-b^-)/(z-b^+)$. Then $T^{-1}(C^+)$ and $T^{-1}(C^-)$
are straight lines which pass through the origin, so $D=T(G)$, where $G$ is a gear
domain (or the image of a gear domain under $z\mapsto1/z$).  For the
converse, when the extended tooth edges intersect in two points, $C^+\cap C^-=\{T(0),T(\infty)\}$ where $T$ is a M\"obius
transformation sending some $G_{\beta,\gamma}$ to the pregear
$D$. \qed

\begin{prop}\label{prop:gearcondition}
  Let $D$ be a pregear.  Then $D$ is a gear if and only if its tooth
  edges are straight, or equivalently, if the A- and B-arcs are
  concentric.
\end{prop}
\subsection{Degenerate gears and pregears}

\begin{defi}\rm
A circular quadrilateral is called a \emph{degenerate} gear (or
pregear) if it is not a gear (or pregear) but is arbitrarily close to
one. 
\end{defi}
Naturally a degenerate gear is also a degenerate pregear. Degenerate
pregears have a very rigidly defined structure. According to
Proposition~\ref{prop:pregearcondition}, each edge of a pregear is a
circular arc orthogonal to the two edges adjacent to it, and the tooth
edges lie in circles $C^\pm$ which intersect in two points.  It
follows from this that in a degenerate pregear the circles $C^\pm$
must be tangent.  Further, any circle orthogonal to both of them must
pass through the point of tangency.  From this one may deduce the
following.

\begin{prop} \label{prop:degeneratepregears} Let $D$ be a bounded
  degenerate pregear symmetric in $\R$.  Then the circles $C^\pm$ containing
  the tooth edges are tangent, say at some point $w^*$, and are orthogonal
  to the A- and B-circles (the circles containing the non-tooth
  edges), which are also tangent at $w^*$.
\end{prop}

\begin{figure}[!htb]
  \centering
 \pic{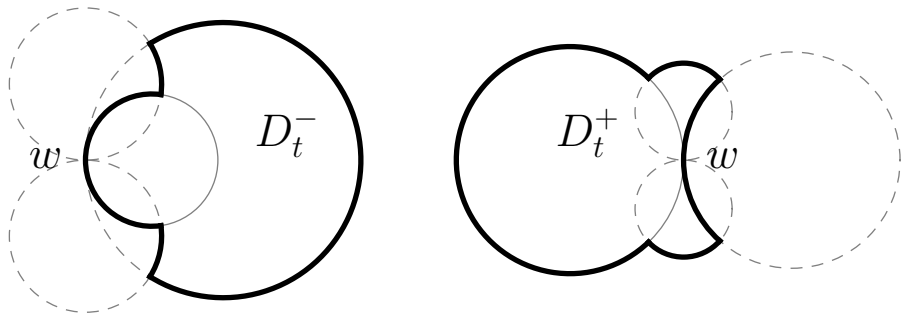}{-.5,0}{3.6}{}
  \caption{The two types of degenerate pregears}
  \label{fig:degeneratepregears}
\end{figure}

The two essentially different possibilities allowed by Proposition
\ref{prop:degeneratepregears} are when the point of tangency $w^*$ of
the tooth edges is at what we have been calling $f(-1)$ or at $f(1)$,
i.e., according to whether the circles containing the A- and B-arcs
are tangent internally or externally, as shown in Figure
\ref{fig:degeneratepregears}.  We will denote by $D_t^-$ and $D_t^+$
respectively the two degenerate pregears with conformal module $M(t)$,
$0<t<\pi$.  For definiteness in the following discussion we normalize
these domains as follows. For $D_t^-$ we apply a real M\"obius
transformation to assure that the A-circle is $\{|w+1/2|=1/2\}$ and
the B-circle is $\{|w|=1\}$, so the point of internal tangency of
$C^\pm$ is at $w=-1\in\partial D_t^-$.  This normalization does not
affect the conformal module.  For $D_t^+$ we will assume that the A-
and B-circles are $\{|w|=1\}$ and $\{|w-2|=1\}$ respectively, with
external tangency at $w=1\in\partial D_t^+$.
 
\begin{figure}[!t]  \centering
  \pic{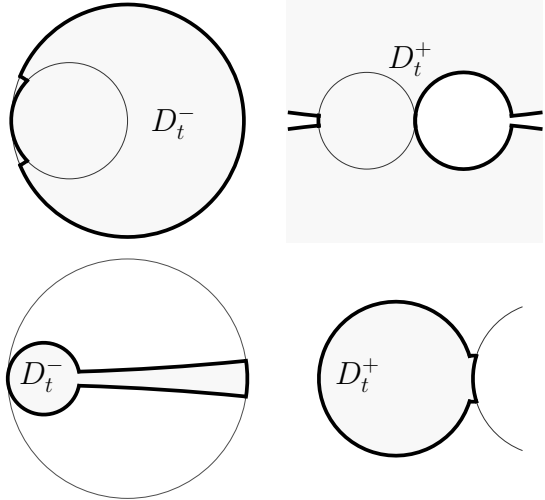}{1,0}{6.3}{scale=.7}
  \caption{Extreme cases of degenerate pregears, for $\eta\to0$ (above) and $\eta\to\infty$ (below).}
  \label{fig:degeneratingdegenerate}
\end{figure}

Let us note the limiting behavior of degenerate pregears as $t\to0$ or
$t\to\pi/2$. Let $\eta$ denote the common radius of $C^\pm$.  Consider
first $D_t^-$. When $\eta\to0$ it can be seen that the conformal
module $M(t)$ tends to $\infty$ (for example, for each $\eta$ apply a
M\"obius transformation leaving the A-circle invariant and sending
$C^\pm$ to circles of radius 1).  In the other direction, when
$\eta\to\infty$, the A-arc is practically all of the A-circle, and the
B-arc is a small arc of the B-circle near $1$, so $M(t)\to0$ (see
Figure \ref{fig:degeneratingdegenerate}).  Similarly, for $D_t^+$, as
$\eta\to0$ we have $M(t)\to\infty$, while as $\eta\to\infty$ we find
that $M(t)\to0$.

\subsection{The boundary of the region of gearlikeness}\label{bounds}

Our main object of study is the the following.
\begin{defi}\emph{
The \textit{region of gearlikeness} is the subset of $\R^2$ defined by 
\begin{eqnarray*}
   \G &=&   \{ (t,\lambda)\colon\ R_{t,\lambda} \mbox{ is the Schwarzian derivative of a
     gear mapping}\}.
\end{eqnarray*}
}
\end{defi}
According to Theorem \ref{theo:lambdainterval},
\[ \G = \{ (t,\lambda)\colon\ 0<t<\frac{\pi}{2},\ 
         \lambda_t^- < \lambda < \lambda_t^+ \}. 
\]
One may obtain very rough bounds on $\lambda_t^-$ and $\lambda_t^+$ by
the classical estimate of Nehari which says that if
$(1-|z|^2)^2|\S_f(z)|>6$ for some $z\in\D$, then the mapping $f$ is
not even univalent.  Applying $z=0$ in (\ref{eq:Rtlambda}) gives the
necessary condition $|16\lambda\cos t + \sin^2 t|<6$ for univalence, from which
\begin{equation}  \label{eq:neharibound}
    \lambda_t^- \ge -\frac{13-\cos 2t}{32\cos t} , \quad
    \lambda_t^+ \le \frac{11+\cos 2t}{32\cos t}.
  \end{equation}  This approach was worked out in
  detail for the analogous case of symmetric quadrilaterals in
  \cite{BrP}, where the Nehari estimate gave rather better results.
  Indeed we give the following exact values for $\lambda_t^+$, $\lambda_t^+$.
 
 \begin{theo} \label{theo:extremlambdas} Let $0<t<\pi/2$.  Then the
  extreme values of $\lambda\in\R$ for which the rational function
  $R_{t,\lambda}$ of (\ref{eq:RtlambdafromRt1t2}) is gearlike are
  given by
  \[ \lambda_t^- = -\frac{1}{4} - \frac{1}{16}\left(\cos t +
    \frac{1}{\cos t}\right) ,\quad
     \lambda_t^+ = \frac{1}{4} - \frac{1}{16}\left(\cos t + 
    \frac{1}{\cos t}\right).
\]  
Further, $R_{t,\lambda_t^-}$, $R_{t,\lambda_t^+}$ are the Schwarzian
derivatives of conformal mappings from the disk\/ $\D$ onto degenerate
pregears of type $D_t^-$, $D_t^+$ respectively, as depicted in Figure
\ref{fig:degeneratepregears}.
\end{theo}

\begin{figure}[!h]
  \centering
  \pic{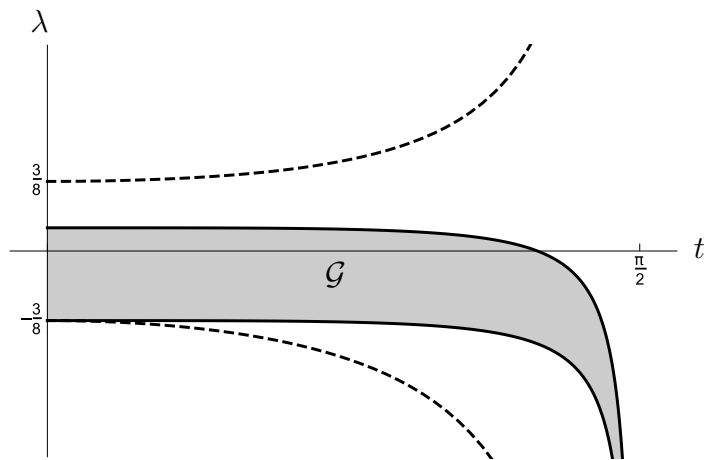}{0.5,0}{5.5}{scale=.7}
  \caption{Region of gearlikeness (gray, bounded by solid curves, which
correspond to degenerate gears).}
  \label{fig:gearlikewithnehari}
\end{figure}

The region of gearlikeness $\G$ is drawn in Figure
\ref{fig:gearlikewithnehari}, together with the rough bound
(\ref{eq:neharibound}).

\vspace{2ex}
\proof For fixed $t$, as $\lambda\to\lambda_t^-$ from above and
$\lambda\to\lambda_t^+$ from below, it is difficult to describe
immediately the behavior of the parameters of the corresponding gear
domains $G_{\beta,\gamma}$.  However, via apropriate M\"obius
transformations we can map these gears to pregears with the A- and
B-circles normalized so that the limiting domains are the degenerate
pregears $D_t^\pm$ (recall Proposition
\ref{prop:degeneratepregears}). Then let us apply a real M\"obius
transformation sending the point of tangency $w^*$ of the extended
edges of the pregear to $\infty$.  Since all image edges are now
straight, the result $\widetilde D_t^\pm$ is an unbounded rectilinear
quadrilateral of a very simple form (Figure
\ref{fig:unboundedquadrilaterals}). The limit mappings
$\D\to\widetilde D_t^\pm$ are given by integrals of
Schwarz-Christoffel type with prevertices at $\pm e^{\pm i t}$, and
also at $-1$ or $1$ which are mapped to vertices of angle $\pi$ at
$\infty$ (cf.\ \cite[Sec.\ 2.1]{DT} for an explanation of this
technical detail):
\begin{figure}[!t]  \centering
  \pic{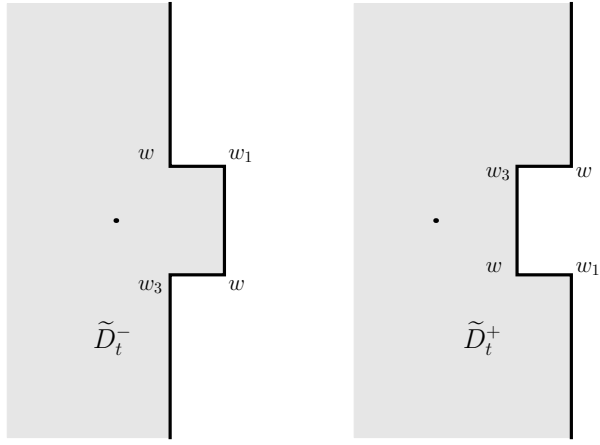}{1,0}{5.5}{scale=.6}
  \caption{Degenerate pregears with tangency mapped to $\infty$}
  \label{fig:unboundedquadrilaterals}
\end{figure}
\begin{eqnarray*}  \label{eq:sch-chr}
 \widetilde f_t^-(z) &=& \int_0^z \frac{1}{(z+1)^2}
   \sqrt{ \frac{z^2 + (2\cos t)z + 1}{z^2 - (2\cos t)z + 1 }  } \,dz,\\
  \widetilde f_t^+(z) &=&  \int_0^z  \frac{1}{(z-1)^2}
   \sqrt{ \frac{z^2 + (2\cos t)z + 1}{z^2 - (2\cos t)z + 1 }  }\,dz.
\end{eqnarray*}
From this the Schwarzian derivatives of the pregear mappings
$f_t^\pm(z)\colon\D\to D_t^\pm$ are equal to those of the
Schwarz-Christoffel integrals $\widetilde f_t^\pm(z)$. These are very easily
calculated, and in the notation of
(\ref{eq:Rtlambda})--(\ref{eq:psi0tpsi1t}) we have
\begin{eqnarray*}  \label{eq:schwdersch-chr}
 \S_{\widetilde f_t^-}(z) &=& -\frac{3+8\cos t+\cos2t}{2(z^4-(2\cos2t)z^2+1)} + \psi_{0,t}(z) ,\\
   \S_{\widetilde f_t^+}(z) &=& -\frac{3-8\cos t+\cos2t}{2(z^4-(2\cos2t)z^2+1)} + \psi_{0,t}(z) 
\end{eqnarray*}
 Upon comparison with
\[  R_{t,\lambda}(z) = \frac{16\lambda\cos t}{z^4-(2\cos2t)z^2+1}+ \psi_{0,t}(z)  
\]
  we obtain the formulas for $\lambda_t^-$ and
$\lambda_t^+$.  
\qed

\section{Discussion and conclusions}

We have described in detail the structure of gear and pregear domains and
shown how the general theory of conformal mapping to circular polygons
can be used to relate the geometry of these domains to the auxiliary
parameters of the corresponding conformal mappings.  In particular, the
conformal module $M(t)$ is a key element for understanding the degeneration
of these domains.

In order to calculate a gear domain numerically by the approach we
have presented here, given a pair of parameters $(t,\lambda)$ it is
necessary to find the appropriate self-mapping $T_q$
of (\ref{eq:Tq}). This amounts to knowing the point $q\in[-1,1]$
which is sent to the gear center.  Since solving the Schwarzian
differential equation (\ref{eq:SfRtlambda0}) a priori produces only a
pregear, it is necessary to find a M\"obius transformation which maps
this to a gear.  This may be approached by solving numerically for the
intersection points of $C^+\cap C^-$ of the circles containing the
tooth edges, or by finding the curvature of the tooth edges
numerically and then adjusting the parameters to assure that this
curvature is zero. Once all this is obtained, the matter of how to
find the unique $(t,\lambda)$ corresponding to prescribed geometric
parameters $(\beta,\gamma)$ can be addressed. These questions will
worked out in the study \cite{BrP3} of numerical aspects of conformal
mappings to gear domains.

\noindent Philip R. Brown\\ 
Department of General Academics\\
Texas A\&M University at Galveston\\ 
PO Box 1675, Galveston, Texas 77553 -1675 \\ 
\texttt{ brownp@tamug.edu }

\medskip
\noindent R. Michael Porter \\
Departamento de Matem\'aticas, CINVESTAV--I.P.N.\\
Apdo.\ Postal 1-798, Arteaga 5 \\
Santiago de Queretaro, Qro., 76000  MEXICO \\
\texttt{mike@math.cinvestav.edu.mx}


\end{document}